\documentclass[a4paper,12pt]{article}

\RequirePackage{amsmath}
\RequirePackage{amssymb}
\def\square{$\vcenter{\hrule\hbox{\vrule height 2truemm \kern 2truemm \vrule }\hrule}$}
\newenvironment{proof}
{\noindent {{\bf Proof}:}}{\hspace*{\fill}\square\vskip 8pt}

\newcommand{\card }{\mathrm{card\,}}

\newcommand{\ord }{\mathrm{ord\,}}

\newcommand{\supp }{\mathrm{supp\,}}

\newcommand{\Div }{{\cal D}\mathrm{iv\,}}
\renewcommand{\div }{\mathrm{div\,}}

\newcommand{\Jac }{\mathrm{ {\cal J}ac}}

\newcommand{\F }{{\mathbb F}}

\newcommand{\ff }{{\mathbb F}}

\newcommand{\ffq }{{\mathbb F}_q}

\newcommand{\fftw }{{\mathbb F}_2}

\renewcommand{\P}{{\mathbb P}}
\newcommand{\A }{{\mathbb A}}

\newcommand{\Z }{{\mathbb Z}}

\newtheorem{prop}{Proposition}
\newtheorem{theo}[prop]{Theorem}
\newtheorem{lem}[prop]{Lemma}
\newtheorem{cor}[prop]{Corollary}

\newtheorem{rem}[prop]{Remark}
\newtheorem{ex}[prop]{Example}

\title{On the Existence of Non-Special Divisors of Degree $g$ and $g-1$ in Algebraic Function Fields over $\F_q$}

\author{St\'ephane Ballet $^{(*)}$ and Dominique Le Brigand $^{(**)}$}
\date{\today}
\begin{document}
\maketitle


\makeatletter


\begin {abstract}
\medskip
 We study the existence of non-special divisors of degree $g$ and $g-1$ for 
algebraic function fields of genus $g\geq 1$ defined over a finite field $\F_q$. 
In particular,
we prove that there always exists an effective non-special  divisor  of degree $g\geq 2$ if  $q\geq 3$
  and
 that
there always exists a non-special  divisor  of degree $g-1\geq 1$ if $q\geq 4$. 
We use our results to
improve upper and upper asymptotic bounds on the bilinear complexity of the multiplication in 
any extension $\F_{q^n}$ of $\F_q$, when $q=2^r\geq 16$.
\end{abstract}

\medskip

{\small \bf 2000 Mathematics Subject Classification:} 11R58.

\renewcommand{\thefootnote}{\fnsymbol{footnote}}
\footnotetext{$^*$ S. Ballet, Laboratoire de G\'{e}om\'{e}trie Alg\'{e}brique 
et Applications \`{a} la
Th\'{e}orie de l'Information,
Universit\'{e} de la Polyn\'{e}sie Fran\c caise, BP 6570, 98702 Faa'a,
Tahiti, Polyn\'{e}sie Fran\c caise,\\
e-mail: ballet@upf.pf.}
\footnotetext{$^{**}$ D. Le Brigand, Institut de Math\'{e}matiques de Jussieu,
Universit\'{e} Pierre et Marie Curie, Paris VI, Campus Chevaleret,
175 rue du Chevaleret, F75013 Paris,\\
e-mail: Dominique.LeBrigand@math.jussieu.fr.}


\section{Introduction}


An important problem in the theory of algebraic function fields is to compute the dimension of a divisor.  In certain cases, it is not an easy task. Moreover, 
given a function field $F/K$ and two integers $n$ and  $d$, 
it is not at all clear if $F$ has a divisor $D$ of degree $d$ with dimension $n$. In fact, the problem occurs when  $0 \leq \deg D \leq 2g_F-2$, where $g_F$ is the genus of $F$. The existence of non-special divisors is, in a sense, related to the number of rational places.
If the full constant field $K$ of an algebraic function field
$F$ is  algebraically closed, then most divisors are non-special and the problem is to find special divisors. Now if $K=\F_q$ is a finite field, the existence of non-special divisors mostly arises for algebraic function fields having few rational places and for $q$ small.
In this paper, we consider an algebraic function field $F/\F_q$ of genus $g$.
We focus on the existence of non-special divisors of degree $d=g$ and 
$d=g-1$ in  $F/\F_q$  because of theorical interest but also because, 
in the case  $d=g-1$, it leads to an improvement of the upper bound  for the bilinear complexity of the multiplication in $\F_{q^n}$ over $\F_q$ when $q=2^r\geq 16$ and also an improvement of the asymptotic bound. 
In fact, this application was our initial motivation.
This paper is organized as follows.
In Section 2, we define the notations and recall basic results.
In Section 3, we give existence results for non-special divisors of degree $g$ and $g-1$.
We settle the problem for $g=1$ and $2$ and 
 we mainly show that if  $q\geq 3$ and $g\geq 2$  (resp. if $q\geq 4$ and $g\geq 2$) there always exist non-special divisors of degree $g$ (resp. $g-1$). 
Finally, in Section 4 we apply the results to the existence of 
 non-special divisors in each step of some towers of function fields. This allows us to
improve upper bounds on the bilinear complexity of the multiplication in any extension of $\F_{2^r}$ when $r\geq 4$.

\section{Preliminaries}

\subsection{Notations}
We mainly use the same notations as in \cite{stich}. Let $F/\F_q$ be an algebraic function field of one variable over $\F_q$.
We assume that the full constant field of $F/\F_q$ is $\F_q$ and denote by $g_F$, or $g$ for short, the genus of $F$. 
Let  $\Div(F/\F_q)$ be the divisor group 
of the algebraic function field $F$ and let $\P(F/\F_q)$, or $\P_F$ for short,  be the set of places of $F$ over $\F_q$. If $u\in F^*$, we denote by $\div(u)$ the principal divisor of $u$ and by $\div_0(u)$ (resp. $\div_\infty(u)$) its zero divisor (resp. pole divisor). Two distinct divisors $D$ and $D'$ are said to be equivalent, denoted $D\sim D'$, if  $D-D'$ is a principal divisor. 
We denote by $\P_k(F/\F_q)$  the set of $k$-degree  places of $F/\F_q$ and by $N_k(F/\F_q)$ (or $N_k$ for short) the order of $\P_k$.
The number $N_1(F/\F_q)$ satisfies the Hasse-Weil inequality 
$ q+1-2g\sqrt{q}\leq N_1(F/\F_q) \leq q+1+2g\sqrt{q}$. 
In particular, if $q$ is a square, $F/\F_q$ is maximal if 
$N_1(F/\F_{q})$ reaches the Hasse-Weil upper bound. 
If ${D}=\sum_{P\in \P_F} n_PP$, we set $\ord_P{D}=n_P$. The support of a divisor ${D}$ is the set 
$\supp(D):=\{P\in \P_F,\, \ord_P{D}\neq 0\}$. 
The divisor ${D}$ is called effective if $\ord_P {D}\geq 0$ for any $P\in \P_F$.
We denote by $\Div^+(F/\F_q)$ the set of effective divisors, by $\A_k$ the set of $k$-degree  effective divisors and set $A_k:=\vert \A_k\vert$. Notice that $A_0=1$ and $A_1=N_1(F/\F_q)$.
We denote by $\Jac(F/\F_q)$ the group of rational points  over $\F_q$ of the jacobian of $F/\F_q$. Then $\Jac(F/\F_q)$ is the group of classes of zero-degree divisors modulo the principal ones and we denote by $[{D}]$ the class of a zero-degree divisor $D$ in $\Jac(F/\F_q)$.
The order $h$ of $\Jac(F/\F_q)$, called divisor class number, is equal to $L(F/\F_q,\,1)$, where $L(F/\F_q,\,t)\in\Z[t]$ is the numerator  of the Zeta function $Z(F/\F_q,\,t)$. Recall that 
$$Z(F/\F_q,\,t):=\sum_{m=0}^{+\infty}A_mt^m= \frac{L(F/\F_q,\,t)}{(1-t)(1-qt)},$$
where $L(F/\F_q,\,t)= \sum_{j=0}^{2g}a_jt^j$,
with $a_j=q^{j-g}a_{2g-j}$, for all $j=0,\ldots , g$.
Let $\pi_j:=\sqrt qe^{i \theta_j}$ and $\bar\pi_j$ be the reciprocal roots of $L$, 
for all $j=1,\ldots ,g$, and
$$L(F/\F_q,\,t)
=\prod_{j=1}^{g}(1-\pi_jt)(1-\bar\pi_jt)
=\prod_{j=1}^{g}(1-2\sqrt q\cos \theta_j\,t+qt^2).$$
The $\theta_j$'s, for $j=1,\ldots ,g$, are the Frobenius angles and we have
\begin{equation}
\label{Lpol}
 \sum_{j=0}^{2g}a_jt^j=\prod_{j=1}^{g}(1-2\sqrt q\cos\theta_j\,t+qt^2).
\end{equation}
Further, we will use the values of the first $a_i$'s:
\begin{equation}
\label{aai}
\left.\begin{array}{rl}
a_1=&N_1-(q+1)\\
a_2=&[N_1^2-(2q+1)N_1]/2+N_2+q\\
a_3=&[N_1^3-3qN_1^2+(3q-1)N_1]/6-(q+1)N_2+N_1N_2+N_3\\
a_4=&[N_1^4+(2-4q)N_1^3-N_1^2-(2-4q)N_1]/24+N_1N_3+N_4+\\
&[(1+2q)N_2+N_2^2
-(1+2q)N_1N_2+N_1^2N_2]/2-(q+1)N_3.\\
\end{array}\right\}
\end{equation}
The real Weil polynomial of $F/\F_q$ is the polynomial $H(T)\in \Z[T] $
defined by:
$$H(T):=\prod_{j=1}^{g}(T-x_j),\mbox{ where } x_j:=-2\sqrt q\cos\theta_j.$$
Using $(\ref{Lpol})$, one can compute the coefficients of $H(T)$ in terms of the $a_j$'s.
Further, we will use the fact that if, for some numerical configuration of the sequence $(N_1,\,N_2,\ldots,\,N_g)$, the corresponding value of $H(2\sqrt q)$ is strictly negative (for instance), then there is no function field of genus $g$ having these numbers of places.\\
If $D\in\Div(F/\F_q)$, then 
$${\cal L}({D}):=\{u \in F^*, \,{D} + \div(u) \geq 0\}\cup \{0\}$$ 
is a $\F_q$-vector space. The dimension of ${\cal L}({D})$, denoted by $\dim D$ and called dimension of divisor $D$, 
is given by Riemann-Roch Theorem 
$$\dim{D}=\deg D-g+1+i(D),$$
where $i({D})$, the index of speciality of ${D}$, is equal to the dimension over $\F_q$ of ${\cal L}({\kappa}- {D})$, $\kappa$ being a canonical divisor. 
 A divisor $D$ is called non-special if $i({D})=0$ and 
otherwise it is called special. 
The index of speciality of a divisor can also be defined in terms of differentials.
The set of regular differentials of $F/\F_q$ is denoted by $\Omega_F(0)$ and one has $\dim_{\F_q}\Omega_F(0)=g$. 

\subsection{Basic results} 
Recall some results about non-special divisors (cf. \cite{stich}).
If $\deg{D}<0$, then  
$\dim{D}=0$ and if $\deg{D}\geq 0$ then $\dim{D} \geq \deg{D}-g+1$.
When   
$0 \leq \deg{D} \leq 2g-2$, the computation of $\dim{D}$ is difficult. Anyway,  one has some general results.

\begin{prop}
\label{basic}
\begin{enumerate}
\item 
$\ffq\subset {\cal L}(D)$ if and only if  $D\geq 0$.
\item
If $\deg{D} > 2g-2$ then ${D}$ is non-special. 
\item
The property of a divisor ${D}$ being special or non-special depends only on  the class of $D$ up to equivalence.
\item
Any canonical divisor $\kappa$ is special, $\deg\kappa=2g-2$  and $\dim\kappa=g$. 
\item
Any divisor ${D}$ with $\dim{D}>0$ and $\deg{D}<g$ is special. 
\item
If ${D}$ is non-special and $D' \geq {D}$, then $D'$ is non-special.  
\item
For any divisor ${D}$ with $0 \leq \deg{D} \leq 2g-2$ holds 
$\dim{D}\leq 1+\frac12\deg{D}.$
\end{enumerate}
\end{prop} 

 \noindent 
 For the rational function field $F=\ffq(x)$ ($g=0$), there is no non-zero regular differential, thus all divisors of degree $d\geq 0$ are non-special. So we assume from now on that $g\geq 1$ and we
 focus on the existence of non-special divisors of degree $g$ or $g-1$.
Note that $g-1$ is the least possible degree for a divisor $D$ to be non-special, since then
$0\leq \dim D= \deg D-g+1$. 
Moreover, if $N_1(F/\F_q)\geq 1$ and if there exists a non-special divisor of degree  $g-1$, then there exists a non-special divisor of any degree $d\geq g-1$ by assertion 6 of Proposition \ref{basic}.
We have the following trivial observations.

\begin{lem} 
\label{trivial}
Assume $g\geq 1$.
Let $D\in \Div(F/\F_q)$ and set $d=\deg D$. 
\begin{enumerate}
\item If $d=g$, $D$ is a non-special divisor if and only if $\dim D=1$. Assume $D$  is a non-special divisor of degree $g$, then $D\sim D_0$, where $D_0$ is effective.
If $D\geq 0$ and $d=g$, 
$D$ is a non-special divisor if and only if ${\cal L}(D)=\ffq$.
\item If $d=g-1$, $D$ is a non-special divisor if and only if $\dim D=0$. A non-special divisor of degree $g-1$, if any, is never effective.
\item If $g>1$ and $A_{g-1}=0$, then any divisor of degree $g-1$  is non-special. 
\end{enumerate}
\end{lem} 

 \noindent
A consequence of assertion 1 of  Lemma \ref{trivial} is:
\begin{lem}
\label{D-P}
 Assume that $D\in \Div(F/\F_q)$ is an effective non-special divisor of degree $g\geq 1$. If there exists a degree one place such that $P\not\in \supp(D)$, then 
$D-P$ is a non-special divisor of degree $g-1$.
\end{lem}

\section{Existence of non-special divisors of degree $g-1$ or $g$}
\subsection{General case}
Let $F/\F _q$ be an algebraic function field of genus $g$.
We denote by ${\cal E}_g$ and ${\cal E}_{g-1}$ the following properties:
$$\begin{array}{ll}
{\cal E}_g\,:&F/\F_q\mbox{ has an effective non-special divisor of degree } g,\\
{\cal E}_{g-1}\,:&F/\F_q\mbox{ has a non-special divisor of degree } g-1.\\
\end{array}$$
If  $F/\F_q$ has enough rational places compared to the genus, then ${\cal E}_g$ and ${\cal E}_{g-1}$ are true.

\begin{prop} 
\label{nonspec1}
Let $F/\F _q$ be an algebraic function field of genus $g\geq 1$.
 \begin{enumerate}
\item If $N_1(F/\F _q)\geq  g$, there exists a non-special divisor 
${D}$ such that $D\geq 0$,  $\deg { D} = g$ and $\supp D\subset \P_1(F/\F _q)$.
\item  If $N_1(F/\F _q)\geq  g+1$, there exists a non-special divisor such that
$\deg {D} = g-1$ and $\supp D\subset \P_1(F/\F _q)$. 
\end{enumerate}
\end{prop}   
\begin{proof}
 \begin{enumerate}
\item cf.  \cite[Proposition I.6.10]{stich}.
\item Let $T\subset \P_1(F/\F _q)$ be such that $\vert T\vert =g$ and,
using assertion 1,  let
$D\geq 0$ be a non-special divisor such that $\deg { D} = g$ and $\supp D\subset T$.
Select $P\in \P_1(F/\F _q)\setminus \supp(D)$ and 
apply Lemma \ref{D-P}.
\end{enumerate}
\end{proof}
 
\noindent

\begin{prop}
\label{$A_{g-1}<h$}
 Let $F/\ffq$ be a function field of genus $g$. 
Denote by $h$  the order of the divisor class group of $F/\ffq$. 
\begin{enumerate}
\item If 
$ A_g<h(q+1)$, then ${\cal E}_{g}$ is true.
\item
If $A_{g-1}<h$, then  
${\cal E}_{g-1}$ is true.
\item Assume $g\geq 2$. 
If $ A_{g-2}<h$, then ${\cal E}_{g}$ is true.
\end{enumerate}
\end{prop}
\begin{proof} Recall that, in any function field, there exists a degree $1$ divisor (this a result of F.K. Schmidt, see \cite[Cor. V.1.11]{stich} for instance), so there exist divisors of any degree $\geq 1$.
Let $d\geq 1$ and $D_0\in \A_d$, if any, and consider the map $\psi_{d,D_0}$:
\begin{equation}
\label{mapD0}
\begin{array}{rcl}
\psi_{d,D_0}: \,\A_{d}& \longrightarrow& \Jac(F/\F_q)\\
{D}&\mapsto&[{D}- {D}_0].\\
\end{array}
\end{equation}
\begin{enumerate}
\item 
First, it is well-known that $1\leq h\leq A_g$ is true for any function field. Indeed, let $D$ be such that $\deg D=g$. By Riemann-Roch, $\dim D\geq 1$, thus there exists an effective divisor of degree $g$ which is equivalent to $D$. So assume $D_0\in\A_g$ and consider the map $\psi_{g,D_0}$.
For all $[R]\in   \Jac(F/\F_q)$, we have $\deg (R+D_0)=g$, thus $\dim (R+D_0)\geq 1$ and there exists $u\in F^*$ such that $D:=R+D_0+\div(u)$ is in $\A_g$ and $[R]=[D-D_0]=\psi_{g,D_0}(D)$. This proves that $\psi_{g,D_0}$ is surjective and that  $h\leq A_g$. Assume now that $F/\F_q$ has no non-special divisor $D$ of degree $g$. Then $\dim D\geq 2$ for all degree $g$ divisors,
thus, for all $[R]\in  \Jac(F/\F_q)$, we have
$$\card\{D\in \A_g,\, [D-D_0]=[R]\}=\frac{q^{\dim(R+D_0)}-1}{q-1}\geq \frac{q^2-1}{q-1}=q+1$$
and $A_g\geq h(q+1)$.
\noindent Notice that this is less restrictive than Lemma 5 of \cite{nichala}, which assumes that $N_1(F/\Fá_q)\geq 1$.

\item 
A divisor  ${D}$   of degree $g-1$  is non-special  if and only if $\dim {D}=0$. 
If $g=1$, there exists a non-special divisor of degree $g-1=0$ if and only if 
$h=N_1>1=A_0$, since two distinct degree one places are not equivalent.
Assume now that $g>1$.
Hence,  it is sufficient to prove the existence of a divisor  of degree $g-1$ 
which is not equivalent to any effective divisor. If $A_{g-1}=0$, the result is proved.
Otherwise, let ${D}_0$  be an effective divisor of degree $g-1\geq 1$ and  consider the map 
$\psi_{g-1,D_0}$.
If  $A_{g-1} < h$, this map is not surjective.  
Hence, there exists a zero-degree divisor $ R$ such that  $[{ R}]$  is not in the image of $\psi_{g-1,D_0}$. Consequently, ${D}={R}+{D}_0$ is a divisor of degree $g-1$ 
which is not equivalent to an effective divisor.   Thus ${D}$ is non-special.

\item From the functional equation of the zeta function, it can be deduced (see \cite[Lemma 3 (i)]{nichala}) that,  for $g\geq 1$, one has
\begin{equation}
\label{zeta1}
A_n=q^{n+1-g}A_{2g-2-n}+h\frac{q^{n+1-g}-1}{q-1},\, \mbox{ for all }0\leq n\leq 2g-2.
\end{equation}
For $g\geq 2$ and $n=g$ this gives
\begin{equation}
\label{Ag}
A_g=h+qA_{g-2}.
\end{equation}
Thus if $g\geq 2$,
$$A_g<(q+1)h\iff A_{g-2}<h.$$

 \end{enumerate}
\end{proof}

\noindent We quote the following consequence of assertion 2.
\begin{cor}
\label{jac}
 Let $F/\F_q$ be an algebraic function field of genus $g\geq 2$ such that $A_{g-1}\geq 1$.
 Denote by $h$  the order of the divisor class group of $F/\ffq$.  
Then ${\cal E}_{g-1}$ is untrue if and only if 
  there exists $h$ elements of $\A_{g-1}$ pairwise non-equivalent.
\end{cor}
\begin{proof}
Let $r$ be the maximum number of pairwise non-equivalent elements of $\A_{g-1}$
and let $D_1,\ldots , D_r$  be elements of $\A_{g-1}$ pairwise non-equivalent. Then
$$\{[0]=[D_1-D_1],\, [D_2-D_1],\ldots,\, [D_r-D_1]\}$$
is a subset of $ \Jac(F/\F_q)$ of order $r$.
If $r=h$, for any divisor $D$ of degree $d=g-1$,  we have $[D-D_1]=[D_i-D_1]$ for  some $i$, $1\leq i\leq h$, and then $D\sim D_i$. Thus $\dim D\geq 1$.
If $r<h$, $\psi_{g-1,D_1}$ is not surjective and the result follows. 

\end{proof}

\subsection{Case $g=1$}\label{g=1}
If the genus of $F/\F_q$ is $g=1$, any divisor of degree $d=g$ 
is non-special since $d\geq 2g-1=1$ and there exists a non-special divisor of degree $g-1=0$ if and only if the divisor class number $h$ is $>1$, i.e. $N_1\geq 2$.
So there are exactly 3 function fields of genus $1$ which have no non-special divisor of degree $g-1$. They are the elliptic solutions to the divisor class number one problem (see \cite{MacR} and  \cite{MQ}): 
$$\begin{array}{crc}
q=2,&y^2+y+(x^3+x+1)=0,&\\
q=3,&y^2-(x^3+2x+2)=0,&\\
q=4,&y^2+y+(x^3+a)=0,&\mbox{ where } \F_4=\F_2(a).\\
\end{array}$$
So, in the rest of this paper, except otherwise stated,
we assume that the genus of a function field is $\geq 2$.

\subsection{Existence of non-special divisors of degree $g\geq 2$}
An algebraic function field of genus $g\geq 2$ has an effective non-special divisor of degree $g$ if it has enough places of degree $2$.

\begin{lem}
\label{N2geq}
 Let $F/\F_q$ be an algebraic function field of genus $g\geq 2$. Then 
${\cal E}_{g}$ is true in either of the following cases
\begin{enumerate}
\item[\rm(i)] 
$N_2\geq q+2$.
\item[\rm(ii)] $N_2=q+1$ and $N_g\geq 1$. 
\item[\rm(iii)] $N_2=q+1$ and $N_1\geq 1$.
\end{enumerate}
\end{lem}
\begin{proof} 
Assume that there is no non-special divisor of degree $g$. Then, by Proposition \ref{$A_{g-1}<h$}, one has $A_{g-2}\geq h$.
But, since  $A_g\geq N_2A_{g-2}$ and $A_g=h+qA_{g-2}$ by
 $(\ref{Ag})$, we have
$$h-A_{g-2}\geq (N_2-q-1)A_{g-2},$$
which contradicts inequality $A_{g-2}\geq h$ 
as soon as $(N_2-q-1)A_{g-2}\geq 1$ and in particular if $N_2\geq q+2$.
We obtain a contradiction  also for $N_2= q+1$ as soon as $A_g>N_2A_{g-2}$, which is the case
if there exists an effective divisor of degree $g$ which support does not contain a degree two place. This is the case if $N_g\geq 1$ or $N_1\geq 1$ or, more generally, if $N_d\geq 1$ where $d\neq 2$ is an integer dividing $g$.
 \end{proof}

\noindent
We want to deduce from Proposition \ref{$A_{g-1}<h$} an existence result for non-special divisors 
of degree $g\geq 2$ which is more general that Lemma 6 of   \cite{nichala}. 
We will use the following result.

\begin{prop}
\label{h=x-pb}
Let $F/\ffq$ be a function field of genus $g\geq 2$. We denote by $h$ its divisor class number.
\begin{itemize}
\item 
Up to isomorphism, there are $4$ function fields $F/\ffq$, $2$ of them being hyperelliptic, 
 such that $h=1$ and $g\geq 2$. They are obtained for $F=\fftw(x,y)$ with
$$\begin{array}{c|l|c|c|c}
g&\mbox{ \rm equation }&N_1&N_2&N_3\\
\hline
2&y^2+y+(x^5+x^3+1)=0&1&2&\\
2 &y^2+y+(x^3+x^2+1)/(x^3+x+1)=0&0&3&\\
3&y^4+xy^3+(x+1)y+(x^4+x^3+1)=0&0&0&1\\
3 &y^4+xy^3+(x+1)y+(x^4+x+1)=0&0&1&1
\end{array}$$
\item
Up to isomorphism, there are $15$ function fields $F/\ffq$, $7$ of them being hyperelliptic, such that $h=2$ and $g\geq 2$.
 They are obtained for $F=\ffq(x,y)$ and\\
$$\begin{array}{c|c|l|c|c|c}
q&g&\mbox{ \rm equation }&N_1&N_2&N_3\\
\hline
3&2&y^2-2(x^2+1)(x^4+2x^3+x+1)=0&0&5&\\
\hline
2&2&y^2+y+(x^3+x+1)/(x^2+x+1)=0&1&3&\\
&&y^2+y+(x^4+x+1)/x=0&2&1&\\
\hline
2&3&y^2+y+(x^4+x^3+x^2+x+1)/(x^3+x+1)=0&1&2&1\\
&& y^2+y+(x^5+x^2+1)/(x^2+x+1)=0&1&3&0\\
&&  y^2+y+(x^6+x+1)/(x^2+x+1)^3=0&0&4&2\\
&& y^2+y+(x^4+x^3+1)/(x^4+x+1)=0&0&3&2\\
&& y^4+xy^3+(x+1)y+(x^4+x^2+1)=0&0&2&2\\
&& y^3+(x^2+x+1)\,y+(x^4+x^3+1)=0&1&0&3\\
&& y^3+y+(x^4+x^3+1)=0&1&1&2\\
\end{array}$$
or $q=2$, $g=4$ and 
$$\begin{array}{l}
y^3+(x^4+x^3+1) y+(x^6+x^3+1)=0  \,,\, (N_j)_{1\leq j\leq g}=(0,0,4,2)\\ 
y^3+(x^4+x^2+1) y+(x^6+x^5+1)=0  \,,\, (N_j)_{1\leq j\leq g}=(0,0,4,2)\\
y^3+(x^4+x^3+1)\, y+(x^6+x+1)=0  \,,\, (N_j)_{1\leq j\leq g}=(0,1,3,3)\\
y^6+xy^5+(x^2+1)y^4+(x^3+x^2)y^3+x^6+x^5+x^3+x+1=0 \,,\, \\
\ \ \ \ \ \ \ \ \ \ \ \ \ \ \ \ \ \ \ \ \ \ \ \ \ \ \ \ \ \ \ \ \ \ \ \ \ \ \ \ \ \ \ \ \ \ \ \ \ \ 
(N_j)_{1\leq j\leq g}=(0,1,1,3)\\
y^6+xy^5+x^3y^3+y^2+(x^5+x^2)y+x^6+x^2+1=0 \,,\, \\
\ \ \ \ \ \ \ \ \ \ \ \ \ \ \ \ \ \ \ \ \ \ \ \ \ \ \ \ \ \ \ \ \ \ \ \ \ \ \ \ \ \ \ \ \ \ \ \ \ \ 
(N_j)_{1\leq j\leq g}=(0,1,2,3)\\
\end{array}$$
\end{itemize}
\end{prop}
\begin{proof}
See \cite{MQ} and  \cite{LMQ}  for  the solutions of the $(h=1)$-problem  and  \cite[Prop. 3.1. and Th. 4.1.]{DLB2}
for  the solutions of the $(h=2)$-problem
\footnote{Note that in
\cite[Th. 4.1 and its proof]{DLB2} there are misprints in the last two equations.}.
\end{proof}

\begin{prop} 
\label{nicha}
An algebraic function field $F/\F_q$ of genus $g\geq 2$ has an effective non-special divisor  
of degree $g$ in the following cases: 
 \begin{enumerate}
\item[\rm(i)] 
$q \geq 3$.
\item[\rm(ii)] $q=2$ and $g= 2$, unless
$F:=\fftw(x,y)$, with\\ $y^2+y+(x^5+x^3+1)=0$ or $y^2+y+(x^3+x^2+1)/(x^3+x+1)$. 
\item[\rm(iii)]  $q=2$ and $g=3$ or $g=4$.
\item[\rm(iv)] $q=2$, $g\geq 5$ and $N_1(F/\F_q)\geq 2 $.
\end{enumerate}
\end{prop}

\begin{proof} 
We set $L(t):=L(F/\F_q,\,t)$.
For $g\geq 2$, it follows from $(\ref{zeta1})$ that (see   \cite{LacDes} or
 \cite[Lemma 3 and proof of Lemma 6]{nichala})
$$\sum_{n= 0}^{g-2}A_nt^n+\sum_{n=  0}^{g-1}q^{g-1-n}A_nt^{2g-2-n}=\frac{L(t)-ht^g}{(1-t)(1-qt)}.$$
Substituting  $t=q^{-1/2}$ in the last identity, we obtain
$$2\sum_{ n= 0}^{g-2}q^{-n/2}A_n+q^{-(g-1)/2}A_{g-1}=\frac{h-q^{g/2}L(q^{-1/2})}
{(q^{1/2}-1)^2q^{(g-1)/2}}$$
and since $L(q^{-1/2})=\prod_{i=1}^g\vert 1-\pi_iq^{-1/2}\vert^2\geq 0$, we have
\begin{equation}
\label{zeta2}
 2\sum_{ n= 0}^{g-2}q^{(g-1-n)/2}A_n+A_{g-1}\leq \frac{h}{(q^{1/2}-1)^2}.
 \end{equation}
Assume that $F/\ffq$ has no non-special divisor of degree $g$. Thus 
by Proposition \ref{$A_{g-1}<h$} and $(\ref{Ag})$, one has 
\begin{equation}
\label{Ag1}
A_{g-2}\geq h.
\end{equation}

\begin{enumerate}

\item $q\geq 3$.
Using $(\ref{zeta2})$, $A_{g-2}\geq h$ implies that
$$2q^{1/2}\leq \frac{1}{(q^{1/2}-1)^2},$$
which is absurd if $q\geq 3$. 

\item $q=2$.
For all $m\geq 1$, one has $A_m\geq N_1A_{m-1}$. Thus if $g\geq 2$
\begin{equation}
\label{Angeq}
A_g\geq N_1A_{g-1}\geq N_1(N_1A_{g-2})=N_1^2A_{g-2}.
\end{equation}
If $q=2$ and $g\geq 2$, 
 the inequality $(\ref{Ag1})$ implies $A_g=h+2A_{g-2}\leq 3A_{g-2}$, thus by $(\ref{Angeq})$
$$N_1^2A_{g-2}\leq A_g\leq 3A_{g-2},$$
which is absurd if $N_1\geq 2$. This proves assertion (iv). We are left with $N_1=0$ or $1$. 

\begin{enumerate}
\item 
If $g=2$,
the inequality $(\ref{Ag1})$ implies that $h=1$, since $A_0=1$. By Proposition \ref{h=x-pb}  there are  only two function fields $F/\F_q$ of genus $2$ such that $h=1$. They are such that $q=2$ and $F=\fftw(x,y)$, with
\begin{enumerate}
\item $y^2+y+(x^5+x^3+1)=0$ and $N_1=1$, $N_2=2$, so $A_2=3$. 
\item $y^2+y+(x^3+x^2+1)/(x^3+x+1)=0$ and $N_1=0$, $N_2=3$, so $A_2=3$. 
\end{enumerate}
Since $h=1$, all divisors of a given degree $d>0$ are equivalent. In particular,
all the divisors of degree $g=2$ are equivalent to 
any divisor of $\A_2$, therefore they are special. 

\item $g=3$. 
By inequality $(\ref{Ag1})$, we have $A_{g-2}=A_1=N_1\geq 1$, thus $N_1=1$ and $h=N_1=1$.
We deduce from Proposition \ref{h=x-pb} that there is no solution.

\item  $g=4$. For instance, consider the five function fields with $g=4$ and $h\leq 2$ in Proposition \ref{h=x-pb}.
They are such that $N_1=0$ and $A_{g-2}=N_2<h=2$, thus they have an effective non-special divisor of degree $g=4$.
It can be verified that all the degree $4$ places in
the first two function fields  are non-special, one degree $4$ place in the third one is non-special and none degree $4$ place is non-special in the last two. More generally, we have
$A_2=\frac{N_1(N_1+1)}2+N_2$ and
$$h=A_4-2A_2=
\left\{\begin{array}{lcl}
N_4+N_3+\frac{N_2^2-N_2}2-1&\mbox{ if }& N_1=1\\
N_4+\frac{N_2^2-3N_2}2&\mbox{ if }& N_1=0.\\
\end{array}\right.$$
We can assume that $h\geq 3$, since the case $h\leq 2$ is settled and then we have $A_2\geq h\geq 3$ by the hypothesis $(\ref{Ag1})$.
If $N_2\geq q+2=4$, ${\cal E}_g$ is true by 
Lemma \ref{N2geq} (i). 
If  $N_2= q+1=3$ and $N_1=0$, then $h=N_4\geq 3$ and ${\cal E}_g$ is true by 
Lemma \ref{N2geq} (ii).  If $N_2= q+1=3$ and $N_1=1$,  ${\cal E}_g$ is true by 
Lemma \ref{N2geq} (iii). 
The last possible case is $N_2= 2$, $N_1=1$ and $A_2=3=h=N_4+N_3$. Let us show that it is impossible. 
The real Weil polynomial of a genus $4$ function field is:
$$H(T)=T^4+a_1T^3+(a_2-4q)T^2+(a_3-3qa_1)T+(a_4-2qa_2+2q^2).$$
Using formulae $(\ref{aai})$ in the case $q=2$, $N_1=1$, $N_2= 2$ and $N_3= 3-N_4$,
we obtain:
$$H(T)=T^4-2T^3-6T^2+(11-N_4)T+3N_4-3.$$
Since $U:=H(2\sqrt 2)=(3-2\sqrt 2)N_4+13-10\sqrt 2$ is strictly negative for all $N_4\leq 3$,
 there is no solution.

\end{enumerate}
\end{enumerate}
\end{proof}

\begin{rem}
\label{g>=5}
We conjecture that a less restrictive result, i.e. without any condition on $N_1(F/\F_q)$, is true if $g\geq 5$ and $q=2$.
\end{rem}

\subsection{Existence of non-special divisors of degree $g-1\geq 1$}
We deduce from  Proposition \ref{$A_{g-1}<h$} an existence result for non-special divisors   of degree $g-1$. 
\begin{theo} 
\label{qgeq 4}
Let $F/\ffq$ be a function field of genus $g\geq 2$. Then ${\cal E}_{g-1}$ is true in the following cases
\begin{enumerate}
\item[\rm(i)]  $q\geq 4$.
\item[\rm(ii)] 
$g= 2$, unless
$F/\ffq:=\fftw(x,y)/\fftw$, with $y^2+y=x^5+x^3+1$ or $y^2+y=(x^4+x+1)/x$.
\item[\rm(iii)]  $q= 2$ or $3$, $g\geq 3$ and $N_1\geq q+1$.
\end{enumerate}
\end{theo}
\begin{proof} Recall that, if $A_{g-1}=0$,  the existence is clear. 

\begin{enumerate}
\item $q\geq 4$.
By $(\ref{zeta2})$, we have for $g\geq 2$
$$A_{g-1}<2q^{(g-1)/2}A_0+A_{g-1}\leq 2\sum_{ n= 0}^{g-2}q^{(g-1-n)/2}A_n+A_{g-1}\leq \frac{h}{(q^{1/2}-1)^2}.$$
Thus, if $q\geq 4$, we have $A_{g-1}<h$ and the result follows from  Proposition \ref{$A_{g-1}<h$}. 

\item $g= 2$.
If $A_{g-1}=N_1<h$, the result follows from Proposition \ref{$A_{g-1}<h$}. This is the case
when $N_1=0$ and then all divisors of degree $g-1$ are non-special.
If $N_1\geq g+1=3$, the result is true by Proposition \ref{nonspec1}. 
The  remaining cases are
$N_1=1$ or $2$ with $h=N_1$. 
\begin{enumerate}
\item $N_1=1$ and $h=1$.  There is a unique function field satisfying these conditions (see Proposition \ref{h=x-pb}). It is $F/\ffq:=\fftw(x,y)/\fftw$, with $y^2+y=x^5+x^3+1$. 
Since $h=1$, all divisors of  degree $g-1=1$ are equivalent to the place of degree $1$, thus they are special.
\item $N_1=2$ and $h=2$. There is a unique function field satisfying these conditions (see Proposition \ref{h=x-pb}). It is $F/\ffq:=\fftw(x,y)/\fftw$, with $y^2+y=(x^4+x+1)/x$.  Since the two degree one places are non-equivalent, it follows from Corollary \ref{jac} that ${\cal E}_{g-1}$ is untrue.
\end{enumerate}

\item  $q=3$, $g\geq 3$.
If $q= 3$ and $g\geq 3$, we have seen in the proof of Proposition \ref{nicha} that 
\begin{equation}
\label{eqqq}
 A_g<(q+1)h.
\end{equation}
Using that $A_g\geq N_1A_{g-1}$,
we obtain $N_1A_{g-1}\leq A_g<4h$
and, if $N_1\geq 4=q+1$, we have $A_{g-1}<h$.
Note that, if $g=3$, it does not improve assertion 2 of Proposition \ref{nonspec1}.

\item $q= 2$, $g\geq 3$.
If $q= 2$, $g\geq 3$ and $N_1\geq 2$, we also proved that 
$(\ref{eqqq})$ is true. Similarly, we obtain $N_1A_{g-1}\leq A_g<3h$,
thus, if $N_1\geq 3=q+1$, we have $A_{g-1}<h$.
\end{enumerate}
We use assertion 2 of Proposition  \ref{$A_{g-1}<h$} to finish the proof.
\end{proof}

\begin{rem} We can  prove that ${\cal E}_{g-1}$ is always true if
$ g=3$ and $q=3$, and, if $g=3$ and $q=2$, there is only a finite number of exceptions, which are the following
$$\begin{array}{l|c|c|c|c}
\mbox{equation}&N_1&N_2&N_3&h\\
\hline
y^4+xy^3+(x+1)y+(x^4+x+1)=0&0&1&1&1\\
 y^2+y+(x^6+x+1)/(x^2+x+1)^3=0&0&4&2&2\\
y^4+xy^3+(x+1)y+(x^4+x^2+1)=0&0&2&2&2\\
y^3+y+(x^4+x^3+1)=0&1&1&2&2\\
y^3+x^2y^2+(x^3+1)y+(x^4+x^3+1)=0&1&2&2&3\\
y^3+x^2y+(x^4+x^3+x)=0&2&0&3&3\\
y^3+(x^2+x+1)y+(x^4+x+1)=0&1&3&2&4\\
\end{array}$$ 

\end{rem}


\subsubsection{Constant field restrictions of maximal function fields}
In the following Lemma, we give the value of $A_{g-1}$
 in terms of the coefficients of the  polynomial  $L(F/\F_q,t)$.
\begin{lem}
\label{A_{g-1}=}
Let $F/\ffq$ be a function field of genus $g$ and let $L(t)= \sum_{i=0}^{2g}a_it^i$ be the numerator of its Zeta function.
Then 
$$A_{g-1}=\frac 1{q-1}\left(h-\left(a_g+2\sum_{i=0}^{g-1}a_i\right)\right ).$$
\end{lem}
\begin{proof}
This is a well-known result (cf. \cite{tsfvl} Section 5 or \cite{pellik}). 
From
$$Z(t)
=\sum_{m=0}^{+\infty}A_mt^m= \frac{L(t)}{(1-t)(1-qt)}
=\frac{ \sum_{i=0}^{2g}a_it^i}{(1-t)(1-qt)}$$
we deduce that for all $m\geq 0$, 
$$A_m=\sum_{i=0}^{m} \frac{q^{m-i+1}-1}{ q-1}a_i.$$
In particular, $$(q-1)A_{g-1}=\sum_{i=0}^{g-1}(q^{g-i}-1)a_i.$$ 
Since 
$a_i=q^{i-g}a_{2g-i}$, for all $i=0,\ldots g$, we get
$$(q-1)A_{g-1}=q^g\sum_{i=0}^{g-1}q^{-i}a_i-\sum_{i=0}^{g-1}a_i=q^g\sum_{i=0}^{g-1}q^{-i}q^{i-g}a_{2g-i}-\sum_{i=0}^{g-1}a_i.$$ 
Hence $$(q-1)A_{g-1}=\sum_{i=0}^{g-1}(a_{2g-i}-a_i).$$
Furthermore, we know that $h=L(1)=\sum_{i=0}^{2g}a_i$, therefore 
$$A_{g-1}=\frac{1 }{ q-1}\left (h-\left (a_g+2\sum_{i=0}^{g-1}a_i\right)\right ).$$ 
\end{proof}

\noindent
Using the preceding Lemma, we obtain a corollary to assertion 2 of 
Proposition \ref{$A_{g-1}<h$}.


\begin{cor}\label{coropi}
If $F/\F_q$ is an algebraic function field such that 
$q\geq 3$ and $a_g+2\sum_{i=0}^{g-1}a_i\geq 0$ 
(resp. $q=2$ and $a_g+2\sum_{i=0}^{g-1}a_i> 0$),
then ${\cal E}_{g-1}$ is true.
\end{cor}
\begin{proof}
By 
Lemma  \ref{A_{g-1}=}, we have $A_{g-1}< h$. The result follows using 
Proposition \ref{$A_{g-1}<h$}.
\end{proof}

\noindent
We will give  examples of function fields satisfying the hypothesis of Corollary \ref{coropi},
but before that, we recall the following result. 

\begin{lem}
\label{constres}
 \begin{enumerate}
\item
Let $F/\F_{q^2}$ be a maximal function field. Then the reciprocal roots of $L(F/\F_{q^2},t)$ are $\pi_i=-q$, for all $i=1,\ldots,2g$, and thus $L(F/\F_{q^2},t)=(1+qt)^{2g}$. 
\item Let $G/\F_{q}$ be a function field such that its constant field extension 
$F/\F_{q^2}=G.\F_{q^2}/\F_{q^2}$ is maximal, then  $L(G/\F_{q},t)=(1+qt^2)^g$.
Moreover, $N_1(G/\F_{q})=q+1$.
\item Reciprocally, if $G/\F_{q}$ is a function field such that $L(G/\F_{q},t)=(1+qt^2)^g$, then its constant field extension 
$F/\F_{q^2}=G.\F_{q^2}/\F_{q^2}$ is maximal.
\end{enumerate}
\end{lem}
\begin{proof}\begin{enumerate}
\item See \cite[Proposition V.3.3.]{stich}.
\item The genera of $G/\F_{q}$ and its constant field extension $F/\F_{q^2}$ are equal.
Let us denote by $\alpha_i$ the reciprocal roots of $L(G/\F_{q},t)$. Then 
the reciprocal roots of  $L(F/\F_{q^2},t)$ are $\pi_i=\alpha_i^2$. Since $\pi_i=-q$, we have $\alpha_i= i\sqrt q$ and $\bar\alpha_i=-i\sqrt q$ and the result follows. In particular,
$N_1(G/\F_{q})=q+1-\sum_{i=1}^g(\alpha_i+\bar\alpha_i)=q+1$.
\item Clear.
\end{enumerate}
\end{proof}

\begin{ex} {\rm The Hermitian function field 
$F/\ff_{q^2}$ is such that $F=\ffq(x,y)$ with $y^q+y-x^{q+1}=0$.
It is a maximal function field of genus $g=\frac{q(q-1)}2$ and it
is the constant field extension of $G/\ffq$, where $G=\ffq(x,y)$,  with $y^q+y-x^{q+1}=0$.
We can say that $G/\ffq$ is a ``constant field restriction" of $F/\ff_{q^2}$.
Lemma \ref{constres} applies to $G/\ffq$. 
Recall that  all subfields $L/\ff_{q^2}$ of the Hermitian function field $F/\ff_{q^2}$ are maximal function fields.
}
\end{ex}

\begin{cor}
\label{maxnspe}
If the algebraic function field $G/\F_q$ is a constant field restriction of a maximal function field
$F/\F_{q^2}=G.\F_{q^2}/\F_{q^2}$, then $G/\F_q$ contains a non-special divisor of degree $g-1$.
\end{cor}
\begin{proof} By the preceding Lemma,  $L(G/\F_{q},t)=(1+qt^2)^g$, thus, by Corollary \ref{coropi}, $G/\F_{q}$  contains a non-special divisor of degree $g-1$.
\end{proof}

\begin{rem}
{\rm 
Corollary  \ref{maxnspe} 
does not improve the previous results. 
Indeed, consider a constant field restriction $G/\ffq$ of a maximal function field 
over $\ff_{q^2}$. Then $N_1(G/\ffq)=q+1$ and $h=(q+1)^g$. Thus,
${\cal E}_{g}$ and ${\cal E}_{g-1}$ are true by
Proposition \ref{nicha} and 
Theorem \ref{qgeq 4} or, if $g=1$,  by Section \ref{g=1}. 
}\end{rem}

\section{Applications} 

\subsection{Case of a Garcia-Stichtenoth  tower } \label{tower}

In this section we study  the 
tower of function fields introduced in \cite{robal}. 
Let us consider the asymptotic good Garcia-Stichtenoth's abelian tower ${\cal T}_1$ over $\F_{q^2}$ (cf. \cite{gast}), 
$${\cal T}_1:=F_{1}\subset F_{2}\subset  F_{3} \subset \ldots$$
such that $F_1/\F_{q^2}:=\F_{q^2}(x_1)/\F_{q^2}$ is the rational function field, 
$F_2/\F_{q^2}$ is the Hermitian function field, and more generally, for all $k\geq 2$, $F_{k+1}/\F_{q^2}$ is defined recursively by
 $$F_{k+1}:=F_{k}(z_{k+1}),$$ 
where $z_{k+1}$ 
satisfies the equation : 
$$z_{k+1}^q+z_{k+1}=x_k^{q+1}\mbox{ with }x_k:=z_k/x_{k-1}\,.$$ 
If $q=p^r$ with $r>1$, we
define the completed 
tower over $\ff_{q^2}$ considered in \cite{ball3}
$${\cal T}_2:=F_{1,0}\subset F_{1,1}\subset \ldots\subset F_{1,r-1} \subset F_{2,0}\subset F_{2,1} \subset\ldots\subset F_{2,r-1}\ldots$$ 
 such that, for all $k\geq 1$,   $F_{k,0}=F_k$, 
$F_{k,s}/F_k$ is a Galois extension of degree $p^s$, for all $s=1,\ldots ,(r-1)$,
 and $[F_{k,s}:F_{k,s-1}]=p$.
 
 \begin{rem} {
 Let $G_k/\ffq$ be defined recursively by $G_1:=\F_{q}(x_1)$ and, for all $k\geq 1$, $G_ {k+1}:=G_{k}(z_{k+1})$, so $G_k$ is the constant field restriction of $F_k$. This is allowed, since all equations are defined over $\ffq$ and the infinite place of the rational function field $G_1:=\F_q(x_1)$ is fully ramified in each step.
 It can be proved that, for all $k\geq 1$ and $s=0,\ldots,\,r-1$, there exists 
 $z_{k+1,s}\in \ffq[z_{k}]$ such that $F_{k+1,s}=F_{k,s}(z_{k+1,s})$.
Thus we can consider the  constant field restriction $G_{k,s}$ of each step $F_{k,s}$ and
the constant field of $G_{k,s}$ is $\ffq$.
 Of course this is quite clear if $r=1$. 
If $r>1$, it is done  in \cite{robal} for $p=2$ and in \cite{robaldlb} for $p$ odd.
}
\end{rem}
If $q=2^r$ with $r>1$, we consider the tower ${\cal T}_3$ over $\F_q$
 studied in \cite{robal}
\begin{equation}
\label{T3}
{\cal T}_3:=G_{1,0} \subset G_{1,1} \subset \ldots\subset G_{1,r-1}\subset G_{2,0}\subset G_{2,1}\subset \ldots\subset G_{2,r-1},\ldots,
\end{equation}
which is related to
the tower ${\cal T}_2$ by
$$F_{k,s}=G_{k,s}\F_{q^2}\,, \quad \mbox{ for all $k\geq 1$ and $s=0,\ldots ,(r-1)$}.$$
Namely $F_{k,s}/\F_{q^2}$ is the constant field extension of $G_{k,s}/\F_q$. 
Notice that  $G_{1,0}/\F_{q}:=\F_{q}(x_1)/\F_{q}$ is the rational function field 
and  $G_{2,0}/\F_{q}$ is the constant restriction of the Hermitian function field.
Each function field 
$G_{1,s}/\F_q$  is the constant restriction  of $F_{1,s}/\F_{q^2}$, which is maximal
since it is a subfield of the Hermitian function field.
Thus the number of rational places of $G_{1,s}/\F_{q}$, for all $s=0,\ldots ,(r-1)$, and  $G_{2,0}/\F_{q}$ equals $q+1$.
Let us denote by $g_{k,s}$ (resp.$N_{k,s}$) the genus (resp. the number of rational places) of the function field $G_{k,s}/\F_q$.

\noindent
Now,  the following result answers  a question of \cite{robal} in a sense which is explained  in Section \ref{compbil}.

\begin{prop} 
\label{nstour}
Assume $q=2^r$.
Then, for any  function field $G_{k,s}/\F_q$ of the tower ${\cal T}_3$, there exists 
a non-special divisor of degree 
$g_{k,s}-1$.  
\end{prop}

\begin{proof}
We have $g_{k,s}\leq g_{l,t}$ for all $1\leq k\leq l$ and $0\leq s\leq t$. Moreover, 
$g_{1,0}=0$, $g_{2,0}=\frac{q(q-1)}2$ and $N_{1,s}=N_{2,0}=q+1$.\\
If $q=2^r\geq 4$, the result follows from Theorem \ref{qgeq 4} (i),
noticing that none of the steps can be the exception of Section \ref{g=1}.\\
 If $q=2$, there is no intermediate step and we set $G_k/\F_2:= G_{k,0}/\F_2$,
 $g_k:=g_{k,0}$ and $n_k:=N_{k,0}$.
 It can be proved, using the results of \cite{gast}, that $n_k\geq 3=q+1$, for all $k\geq 1$.
 In fact, $n_1=n_2=q+1=3$. The places of $G_2/\F_{2}$ are the pole of $x$, 
 the common zero of $x,z_2$ and the common zero of $x,z_2+1$.
 For all $k\geq 3$, the common zero of $x,z_2,\ldots, z_{k-1}$ in $G_{k-1}/\F_2$
 splits in 
 $G_k/\F_2$ giving two degree one places,
 one being the common zero of $x,z_2,\ldots, z_{k-1},z_k$ and the other being the common zero of $x,z_2,\ldots, z_{k-1},z_k+1$. With the pole of $x$,  
 we obtained at least three degree one places in $G_k/\F_2$.
 Once again, none of the steps can be the exception of Section \ref{g=1} ($g=1$) or the  exception of Theorem \ref{qgeq 4} ($g=2$) and then the result follows.
\end{proof}


\subsection{Previously known applications}
We quote previous works in which  non-special divisors are needed or constructed. 
We must say that,  the existence of such divisors is often clear because the function fields, which are involved,  have many rational places but the problem lies in their effective determination.
\begin{enumerate}
\item In the construction of Goppa codes $C_L(G,D)$ on a function field $F/\ffq$, where $D:=P_1+\cdots+P_n$ is a sum of $n$ distinct places of degree one and $G$ is a divisor, such that $\supp G\cap \{P_1,\ldots,P_n\}=\emptyset$, 
it is often assumed that ${\cal L}(G-D)=\{0\}$. In fact, $C_L(G,D)$ is the image of ${\cal L}(G)$ by the evaluation map
$$\begin{array}{rcl}
 {\cal L}(G)&\rightarrow&\ffq^n\\
 u&\mapsto&(u(P_1),\ldots,u(P_n)).
 \end{array}$$
and this map is injective thanks to the condition ${\cal L}(G-D)=\{0\}$. Thus the dimension of $C_L(G,D)$ is $\dim G$. Of course, any $G$ with  $\deg G<n:=\deg D$ satisfies the condition. But if one wants to consider higher degree, it may be useful to know that there exists a non-special divisor $B$ of degree $g_F-1$, since then, divisor $G':=B+D\geq B$ is also non-special and thus we know the value of $\dim G'$. 
Moreover, there exists $G\sim G'$ such that $\supp G\cap \{P_1,\ldots,P_n\}=\emptyset$ and ${\cal L}(G-D)={\cal L}(B)=\{0\}$.

\item In many constructions of algebraic-geometry codes  (see  \cite{NX} for instance)
or in construction of $(t,s)$-sequences using function fields,  
the existence of a non-special divisor of degree $g$ is assumed. In both cases,  the following basic argument is more or less needed.
Let  $F/\ffq$ be a function field, let
$D:=P_1+\cdots+P_n$ be a sum of $n$ distinct places of degree one and let $G$ be an effective non-special divisor of degree $g_F$. Then, there exists a function $f_i\in {\cal L}(G+P_i)\setminus {\cal L}(G)$ for all $1\leq i\leq n$
and  $(1,f_1,\ldots,f_n)$ is a basis of ${\cal L}(G+D)$.
\item
In \cite{GS2} another asymptotic good tower of function fields is given, which is a sub-tower of ${\cal T}_1$, 
and in \cite{PST} an explicit non-special divisor of degree $g$  is given for each steps. Let us recall the situation. Consider the tower 
 ${\cal F}$ over $\F_{q^2}$ 
$${\cal F}:=F_{1}\subset F_{2}\subset  F_{3} \subset \ldots$$
such that $F_n:=\F_{q^2}(x_1,\ldots,x_n)$,
with  
$$x_{k+1}^q+x_{k+1}=\frac{x_i^q}{x_i^{q-1}+1}\mbox{ for }k=1,\ldots n-1\,.$$ 
Then $[F_n:F_1]=q^{n-1}$ and the infinite place  of the rational function field $F_1/\F_{q^2}:=\F_{q^2}(x_1)/\F_{q^2}$ is fully ramified in each $F_n$. We denote by $P_\infty^{(n)}$ the corresponding place in $F_n$. In  \cite{PST} and for each $n\geq 1$, the authors define explicitly an effective divisor $A^{(n)}$ of $F_n/\F_{q^2}$, such that 
$\dim (c_nP_\infty^{(n)}-A^{(n)})=1$, with $c_n:=g_{F_n}+\deg A^{(n)}$. So 
$D^{(n)}:=c_nP_\infty^{(n)}-A^{(n)}$ is a non-special divisor of degree $g_{F_n}$. 
More precisely, a basis of ${\cal L}(D^{(n)})$ is $\pi_n$ such that
$\div_\infty(\pi_n)=c_nP_\infty^{(n)}$.
We remark that it is then straightforward to show that $D_1^{(n)}:=(c_n-1)P_\infty^{(n)}-A^{(n)}$ is a non-special divisor of degree $g_{F_n}-1$, for each $n\geq 1$,
since ${\cal L}(D_1^{(n)})\subset {\cal L}(D^{(n)})$ and $\pi_n\not\in {\cal L}(D_1^{(n)})$.
\item
In \cite{xing}, the existence of certain divisors $D$ is necessary to obtain asymptotic bounds on frameproof codes. The bound is  obtained using asymptotic good towers of function fields.
We will not recall the definition of frameproof codes. Given a function field $F/\ffq$ of genus $g$ and $P_1,\ldots,P_n$ distinct degree-one places,
 the author of \cite{xing} assumes the existence of an effective divisor $D$ such that $m:=\deg D$
 and $\dim(sD-\sum_{i=1}^nP_i)=0$, for an integer $s\geq 2$. Clearly the greatest possible value for $sm-n$ is then equal to
$g_F-1$ and we note that, if $sm-n=g_F-1$, the divisor $(sD-\sum_{i=1}^nP_i)$ is  
non-special of degree $g_F-1$. However, the existence of such a divisor is hard to prove and the author of \cite{xing} gives a sufficient condition to obtain the result.
\end{enumerate}
We end this (non-exhaustive) enumeration and come to the initial motivation of this work.

\subsection{Application to the bilinear complexity of multiplication}
\label{compbil}
\subsubsection{Problem}
Let $\F_q$ be a finite field with $q$ elements where $q$ is a prime power 
and let $\F_{q^n}$ be a degree $n$ extension of $\F_q$. 
We denote by $m$ the ordinary multiplication  in the finite field $\F_{q^n}$. 
This field will be considered as a $\F_q$-vector space. 
The multiplication $m$ is a bilinear map from $\F_{q^n} \times \F_{q^n}$ into $\F_{q^n}$, 
thus it corresponds to a linear map $M$ from the tensor product 
$\F_{q^n} \bigotimes \F_{q^n}$ over $\F_q$ into $\F_{q^n}$. 
One can also represent 
$M$ by a tensor $t_M \in \F_{q^n}^*\bigotimes \F_{q^n}^* \bigotimes \F_{q^n}$ 
where $\F_{q^n}^*$ denotes the dual of $\F_{q^n}$ over $\F_q$. 
Hence the product of two elements $x$ and $y$ of $\F_{q^n}$ 
is the convolution of this tensor with $x \otimes y \in \F_{q^n} \bigotimes \F_{q^n}$. 
If 
\begin{eqnarray}
\label{tensor}
t_M = \sum^{\lambda}_{l=1}a_l \otimes b_l \otimes c_l,
\end{eqnarray}
where $a_l \in  \F_{q^n}^*$, $b_l \in  \F_{q^n}^*$, $c_l \in  \F_{q^n}$, then 
\begin{eqnarray}
\label{bilmult}
x.y=\sum^{\lambda}_{l=1}a_l(x)b_l(y)c_l.
\end{eqnarray}
 
\bigskip
\noindent
Every expression $(\ref{bilmult})$ is called a bilinear multiplication algorithm $\cal U$. 
The integer $\lambda$ is called the multiplicative complexity  of $\cal U$ and denoted by  $\mu({\cal U})$.

\bigskip
\noindent
Let us set $$ \mu_{q}(n)= \min_{\cal U} \mu({\cal U}), $$  where $\cal U$ 
is running over all bilinear multiplication algorithms in $\F_{q^n}$ over $\F_q$. 
Then $\mu_{q}(n)$ is called the bilinear complexity of 
multiplication in $\F_{q^n}$ over $\F_q$, and it corresponds 
to the least
possible number of summands in any tensor decomposition of type $(\ref{tensor})$.   

\medskip

\subsubsection{An improvement of a multiplication bilinear complexity bound}

\noindent
We know by \cite[Th. 4.2.]{robal} that for $p=2$ and
$q=p^r \geq 16$, the bilinear complexity $\mu _q(n)$
of multiplication in any finite field $\F_{q^n}$ over $\F_q$ satisfies:
 
$$\mu _q(n) \leq 3\left(1+\frac{4p}{q-5}\right) n.$$

\noindent
The above result is obtained in \cite{robal} using the tower ${\cal T}_3$ defined by $(\ref{T3})$.
As  said in  \cite[Section 5]{robal}, the existence of non-special divisors of degree $g_{k,s}-1$ 
for each step $G_{k,s}$ of the tower ${\cal T}_3$, which is proved in  
Proposition \ref{nstour}, enables us to obtain a better bound 
and a better asymptotic bound than the ones obtained in \cite{robal}. 
More precisely, 
we obtain the following Theorem:
\begin{theo}\label{BetterBound}
Assume $q=2^r\geq 16$ and let $n\geq 1$ be an integer.
The bilinear complexity $\mu_q(n)$ of the multiplication
in $\F_{q^n}$ over $\F_q$ satifies
$$ \mu_q(n)\leq 3\left(1+\frac{4}{q-3}\right)n.$$
Thus
$${\cal M}_q=\limsup_{n \rightarrow \infty} \frac{\mu_q(n)}n \leq 3\left(1+
 \frac{4}{q-3}\right).$$
\end{theo}





\begin{rem}
Let us remark that this improvement is obtained using Theorem \ref{qgeq 4} and cannot  
be easily obtained from \cite[Lemma 6]{nichala} as was suggested in \cite[Question 2]{robal}.
Moreover, note that there is a misprint in the upper bound of $\mu_q(n)$ in \cite[Question 2]{robal}. 
In fact, for $p=2$ and $q=2^r\geq 16$, the inequality must be $\mu_q(n)\leq 3\left(1+\frac{2p}{q-3}\right)n.$
\end{rem}

\end{document}